\documentclass{article}
\usepackage{amssymb}
\usepackage{latexsym}
\newtheorem{thm}{Theorem}
\newtheorem{lem}{Lemma}
\newtheorem{cor}{Corollary}
\begin{document}
\title{Special Lagrangian submanifolds and Algebraic Complexity one Torus
Actions}
\author{Edward Goldstein}
\maketitle
\renewcommand{\abstractname}{Abstract}
\begin{abstract}
In the first part of this paper we consider compact algebraic manifolds $M^{2n}$ with an algebraic $(n-1)$-Torus action. We show that there is a $T$-invariant
meromorphic section $\sigma$ of the canonical bundle of $M$. Any such $\sigma$ defines a divisor $D$. On the complement $M'=M-D$ we have a trivialization of the canonical bundle and a $T$-action. If $H^1(M',\mathbb{R})=0$ then results of \cite{Gold2} show that there is a Special Lagrangian (SLag) fibration on $M'$. We will study how the fibers compactify in $M$ and give examples of SLag fibrations on $M'$, including some cases there $H^1(M',\mathbb{R}) \neq 0$. 

In the second part of the paper we study Calabi-Yau hypersurfaces in $M$. We will assume that $\sigma$ is an inverse of a holomorphic section $\eta$ of the anti-canonical bundle of $M$. We will see that under certain assumptions one can choose a holomorphic volume form on the smooth part $D'$ of the divisor $D$  s.t. orbits of the $T$-action give a SLag fibration on $D'$ with respect to the metric, induced from $M$. Transversal sections $\eta_j$ near $\eta$ define smooth Calabi-Yau hypersurfaces $D_j$ in $M$. We will show that one can deform the SLag fibration on $D'$ to SLag fibrations on large parts of $D_j$. This construction applies for instance for $D_j$ being quintics in $\mathbb{C}P^4$ or Calabi-Yau hypersurfaces in the Grassmanian $G(2,4)$.
\end{abstract}
\section{Introduction}
In our previous paper \cite{Gold2} we have shown how to use torus actions
on non-compact Calabi-Yau manifolds to construct SLag
submanifolds. In all the ensuing discussions by SLag submanifolds we mean the following: Let $(M^{2n},\omega,\varphi)$ be a complex manifold with a Kahler form $\omega$ and a non-vanishing holomorphic $(n,0)$-form $\varphi$. Then a submanifold $L^n$ of $M$ is SLag if it satisfies $\omega|_{L}=0 ~ , ~ Im\varphi|_{L}=0$. We refer the reader to \cite{Gold1} and \cite{Mc} for a discussion of SLag submanifolds. In the first part of this paper we develop one general setup, which generalizes all examples in \cite{Gold2}:

Consider a compact algebraic manifold $M^{2n}$ with an algebraic $(n-1)$-torus
action (this is an algebraic complexity one space). We will show (see Theorem 1) that
there is a meromorphic section $\sigma$ of the canonical bundle of $M$,
which is invariant under the $T$-action. Any such $\sigma$ defines a divisor $D$ and on the complement $M'=M-D$ $\sigma$ gives a trivialization of the canonical
bundle. Let $H^1(M',\mathbb{R})=0$. Take any $T$-invariant Kahler form
$\omega$ on $M'$. From \cite {Gold2} we deduce that
$M'$ has a SLag fibration with respect to $\omega$. We will investigate
how the
fibers compactify in $M$. We treat 2 cases : If $\omega$ is the
induced metric from $M$, then Lemma 1 of Section 2.1 covers ``most'' fibers; For $\omega$ a Ricci-flat metric on $M'$ we study one family of examples in Section 2.2.

We will also give examples of SLag fibrations on $M'$, including some cases there $H^1(M',\mathbb{R}) \neq 0$.

In the second part of this paper we investigate Calabi-Yau hypersurfaces in $M$. We assume that the anti-canonical bundle of $M$ is ample and has a $T$-invariant holomorphic section $\eta$ (with a corresponding meromorphic section $\sigma=\eta^{-1}$ of the canonical bundle). We also assume that near the smooth part $D'$ of the zero set $D$ $\eta$ is transversal to zero. This gives a natural trivialization $\varphi'$ of the canonical bundle of $D'$. Consider the induced metric $\omega$ from $M$ on $D'$. We assume that the $T$-action on each connected component $D_0$ of $D'$ has a finite stabilizer. For $D_0$
 one can choose an angle $\theta_0$ s.t.
orbits of the $T$-action are SLag submanifolds of $(D_0,\omega,e^{i\theta_0} \varphi')$. Consider transversal sections $\eta_j$, which converge to $\eta$. Their zero sets $D_j$ are smooth Calabi-Yau hypersurfaces. Near compact, $T$-invariant subsets $N$ of $D_0$ the manifolds $D_j$ converge to $N$ as $\eta_j \mapsto \eta$. We will show (see Theorem 3) that one can perturb our SLag fibration on $N$ to SLag fibrations on the corresponding subsets of $D_j$ (all fibrations are SLag for metrics, induced from $M$). 

We will see that this construction applies for $M=\mathbb{C}P^4$ and $D_j$ being quintics; and for $M=G(2,4)$ (the Grassmanian of 2-planes in $\mathbb{C}^4$) and $D_j$ being Calabi-Yau hypersurfaces in $M$.                          

We will also show that this can generalized from hypersurfaces to complete intersection Calabi-Yau submanifolds of a Kahler manifold $M$ with a higher complexity torus action (see Corollary 1 and Section 3.4).

{\bf Acknowledgments} : This paper is a part of author's work towards his
Ph.D. at MIT. The author is grateful to
 his advisor, Tom Mrowka, for continuing support; and to Shing-Tung Yau for a helpful discussion.

Research is partially supported by an NSERC PGS B Award
\section {Algebraic Complexity one Spaces}
\subsection{Fibration on M' and it's compactification} 
Let $M^{2n}$ be a compact algebraic manifold with an effective algebraic $T^s$-action. The fact the the action is algebraic is equivalent to saying that there is a very ample line bundle $L$ on $M$, and $T$ acts on it's total space. We can use $L$ to construct a meromorphic section of the canonical bundle $K$ of $M$:
\begin{thm}
For $k$ large enough there is a holomorphic section $\sigma_1$ of $K \otimes L^{\otimes k}$ and $\sigma_2$ of $L^{\otimes k}$ s.t. the section $\sigma= \sigma_1 \otimes \sigma_2^{-1}$ is a meromorphic, $T$-invariant section of $K$.
\end{thm} 
{\bf Proof:} Consider the representation $\rho$ of $T$ on $H^0(L,M)$.
 We split it into a direct sum of irreducible representations, i.e. we find sections $\alpha_1,\ldots,\alpha_d$ of $L$ s.t. $T$ acts on $span(\alpha_i)$ with a character $\xi_i$. W.l.o.g we can assume that $\xi_1$ is a trivial character (otherwise we divide the action of $T$ on $L$ by $\xi_1$).

We can view $\xi_j$ as living in the character lattice of the dual Lie algebra
${\cal G}^{\ast}$ of $T$. The condition that the $T$-action is effective on $M$ implies the $\xi_i$ $\mathbb{Z}$-span the character lattice of ${\cal G}^{\ast}$. Indeed suppose $\xi_i$ do not span the character lattice of ${\cal G}^{\ast}$. Then we can find an element $e \neq 1$ of $T$ s.t. $\xi_i(e)=1$, i.e. $e$ acts trivially on $H^0(L,M)$. But then $e$ acts trivially on $M$. Indeed suppose $m \in M$ and $e(m)=m' \neq m$. Since $L$ is very ample, there is a section $\beta$ of $L$ s.t. $\beta(m)=0 ~ , ~ \beta(m') \neq 0$. But then $e$ doesn't act trivially on $\beta$- a contradiction.

For $l$ large enough  $L^{\otimes l} \otimes K$ has a holomorphic section. So $H^0(L^{\otimes l} \otimes K,M)$ is nontrivial. We can find a section $\lambda$ of $L^{\otimes l} \otimes K$ s.t. $V=span(\lambda)$ is $T$-invariant and $T$ acts on $V$ with a character $\phi$. We can write $\phi = \Sigma a_i \xi_i ~ , ~ a_i \in \mathbb{Z}$. We will consider only $a_i \neq 0$ and we divide them into subsets $a_i^+$ and $a_i^-$ of positive and negative values. We will have a section $\sigma_1'=\lambda \otimes \alpha_i^{\otimes -a_i^-}$ of $L^{\otimes k_1} \otimes K$ and a section $\sigma_2'= \otimes \alpha_i^{\otimes a_i^+}$ of $L^{\otimes k_2}$. The $T$-actions leave $V_i=span(\sigma_i')$ invariant and $T$ acts on $V_1$ and $V_2$ with the same character. Also $T$ preserves $\alpha_1$. So tensoring $\sigma_1'$ and $\sigma_2'$ with powers of $\alpha_1$ we deduce that for $k$ large enough we can find sections $\sigma_1$ of $L^{\otimes k} \otimes K$ and $\sigma_2$ of $L^{\otimes k}$ s.t. $T$ acts on $\sigma_i$ with the same character.           Q.E.D. 

From now on we assume that $s=n-1$, i.e. $M$ is an algebraic complexity one space. We refer the reader to \cite{KT} and \cite{Tim} for discussion of complexity one spaces.
Let $\sigma$ be any $T$-invariant meromorphic section of the canonical bundle of $M$. It defines a divisor $D$. Let $M'=M-D$. Take $\omega$ to be any $T$-invariant Kahler form on $M'$ s.t. the $T$-action is Hamiltonian. We assume that $H^1(M',\mathbb{R})=0$. From \cite{Gold2} we get the following result (a similar result was also obtained in \cite{MG2}):
\begin{thm} (\cite{Gold2}) $M'$ has a calibrated fibration $\alpha$ over an open subset of $\mathbb{R}^n$. Moreover for a generic point $p$ (outside of a countable union of $(n-2)$-planes in $\mathbb{R}^n$) the fiber $\alpha^{-1}(p)$ is a smooth SLag submanifold. Connected components of smooth fibers are diffeomorphic to $T^{n-1} \times \mathbb{R}$. If $m \in M'$ is
a regular point of the torus action (i.e. the differential of the action
is injective at $m$), then $m$ is a smooth point of the fiber $L_m$
through $m$ and the fiber $L_m$ is a SLag submanifold near $m$. Singular fibers have singularities of codimension at least 2 and near singular points they are diffeomorphic to a product of a cone with a Euclidean ball.
\end{thm}
The fibration is defined as follows: Let $e_1,\ldots,e_{n-1}$ be a basis for the Lie Algebra of $T^{n-1}$ and $X_1,\ldots,X_{n-1}$ be the corresponding flow vector fields on $M$. Let $\sigma'=i_{X_1} \ldots i_{X_{n-1}} \sigma$ be the $(1,0)$-form on $M'$, obtained by contraction of $\sigma$ by vector fields $X_1,\ldots,X_{n-1}$. Then from \cite{Gold2} we know that $\sigma'$ is a holomorphic $(1,0)$-form on $M'$ and $\sigma=df$ for a $T$-invariant, holomorphic function $f$ on $M'$. Let $\mu$ be a moment map for the $T$-action on $M'$. Then the fibers of the SLag fibration are given as level sets of the function $\alpha=(\mu,Imf)$.

The closure of a connected component $L'$ of any smooth SLag fiber $L$ on $M'$ must intersect $D$ (for otherwise $L'$ would be compact, but it is diffeomorphic to $T^{n-1} \times \mathbb{R}$). 
In order to further understand how the fibers of our SLag fibration on $M'$ compactify in $M$ we will have to make assumptions on the Kahler form $\omega$. First we treat the case then $\omega$ is induced from a Kahler form on $M$.

We decompose $D$ into $D=D_+ \bigcup D_-$
corresponding to the meromorphic and holomorphic parts of $\sigma$. Consider the function $f$ near a point $d$ in $D-D_-$. Then the differential $df=\sigma$ is bounded, and so will be $f$. By Riemann extension theorem, $f$ extends to a holomorphic function near $m$. So $f$ extends to a holomorphic function on $M-D_-$.
 Also $\mu$ is smooth on $M-D_-$ and so our fibration extends to a fibration by level sets of $(\mu,Imf)$ on $M-D_-$. Finally we wish to understand how the fibers compactify near $D_-$. We have the following Lemma:
\begin{lem}
Let $\nu$ be a value of the moment map s.t. $T$ acts freely on $D_{\nu}=\mu^{-1}(\nu) \bigcap D_-$. Consider a fiber $L_t=(\mu=\nu,Imf=t)$ of the fibration over $M-D_-$ for some $t \in \mathbb{R}$. Then the boundary of $L_t$ in $M$ is $D_{\nu}$
\end{lem} 
{\bf Proof}: Obviously the boundary of $L_t$ is contained in $D_{\nu}$. To prove the other inclusion we first note that the $T^{n-1}$ action on $M$ induces the complex torus $T^c=(\mathbb{C}^{\ast})^{n-1}$ action on $M$. Indeed the flow vector fields $X_v$ give rise to the vector fields $JX_v$. The flow of $X_u$ preserves $X_v$ and commutes with $J$, hence it preserves $JX_v$, i.e. $[X_u,JX_v]=0$. Finally the vector field $[X_u -iJX_u,X_v-iJX_v]$ is of type $(1,0)$, but it's equal to $-[JX_u,JX_v]$. So $[JX_u,JX_v]=0$ and all the vector fields $X_u, JX_v$ commute, and this induces the $T^c$-action on $M$ (with the real torus $T^{n-1}$ being the product $S^1 \times \ldots \times S^1$ of unit circles in $T^c$). We have \[{\cal L}_{JX_v}\sigma=d(i_{JX_v}\sigma)=id(i_{X_v}\sigma)=i{\cal L}_{X_v}\sigma=0 \]
So the $T^c$-action preserves $\sigma$, and so it preserves $D,D_+,D_-$. Also since $f$ is holomorphic then $JX_v(f)=iX_v(f)=0$, i.e. $T^c$ preserves $f$. 

Let now $d \in D_{\nu}$. Since $T$ acts freely on $D_{\nu}$ then in particular the differential of the $T$-action on $d$ is injective. So the orbit $L^c$ of $T^c$-action on $d$ is a smooth complex submanifold of $D_-$ of complex dimension $n-1$.  
Since $T^c$-action preserves $D_-$ then it is clear that near $d$ $D_-$ coincides with $L^c$. Also the differential of $\mu|_{L_c}$ is surjective, so near $d$ $D_{\nu}$ coincides with the orbit $L$ of the $T^{n-1}$-action on $d$. 

Consider the level set $\Sigma_{\nu}$ of $\mu$ on $M$. Then near $D_{\nu}$ $\Sigma_{\nu}$ is smooth. Also $T^{n-1}$ acts freely on $D_{\nu}$, hence it acts freely on $\Sigma_{\nu}$ near $D_{\nu}$. We consider the symplectic reduction $M_{red}$ of $\Sigma_{\nu}$ in a $T$-invariant neighbourhood $U$ of $D_{\nu}$. Then $M_{red}$ will be a smooth complex 1-dimensional manifold. Let $\pi: \Sigma_{\nu} \mapsto M_{red}$ be the quotient map. Then $\pi(L)=d'$ is one point in $M_{red}$ (here $L$ is the $T^{n-1}$-orbit of $d$). Also in a sufficiently small neighbourhood $U'$ of $d'$ in $M_{red}$, $\pi^{-1}(U'-d') \bigcap D_- = \emptyset$.

The function $f$ descends to a holomorphic function on $U'$. Suppose that $f$ has a singularity in $d'$. Then the image of $f$ covers a neighbourhood of $\infty$ in $\mathbb{C}$. Hence the image of $Imf$ assumes all real values on $U'-d$, in particular it assumes the value $t$. Hence the SLag fiber $L_t$ in $M-D_-$ intersects the neighbourhood $\pi^{-1}(U')$. From this we easily deduce that $d$ is in the closure of $L_t$ in $M$.

So we need to prove that $f$ has a singularity at $d'$. Suppose not. Then $|f|$ is bounded by some constant $C$. So $|f|<C$ on $\Sigma_{\nu}-D_-$ near $d$. Now the $T^c$-action preserves $f$. The orbit of $\Sigma_{\nu}$ under the $T^c$-action fills a neighbourhood $W$ of $d$ in $M$. On $W-D_-$ we have $|f|<C$. By Riemann extension theorem $f$ extends to a holomorphic function on $W$. So the norm of it's differential $|df|<C'$ for some other constant $C'$, i.e. $|\sigma'|<C'$. But $\sigma'=i_{X_1} \ldots i_{X_{n-1}}\sigma$. We choose a local section $\varphi$ of the canonical bundle of $M$ near $d$ s.t. $|\varphi|=1$. Then $\sigma=g\varphi$ and $|g|$ is not bounded since $\sigma$ has a singularity at $d$. But the vector fields $X_i$ are linearly independent and $\omega$-orthogonal to each other. Hence one easily deduces that $|i_{X_1}\ldots i_{X_{n-1}}\varphi|$ is uniformly bounded from below - a contradiction. So $f$ is singular at $d'$ and we are done.   Q.E.D. 

{\bf Remark 1}: The conclusions of Lemma 1 apply only for those values of $\mu$ on $D_-$, on whose level sets $T^{n-1}$ acts freely. In example 1 in Section 2.2 we'll have that for all non-regular values $\nu$ on $D_-$ the level set $D_{\nu}$ will not intersect the closure of any SLag fiber in $M-D_-$.

{\bf Ricci-Flat metrics} Another interesting class of metrics on $M'$ are $T$-invariant, complete Ricci-flat metrics. We refer the reader to \cite{TY1}, \cite{Joy} for some existence results. If $H^1(M',\mathbb{R})=0$ then from Theorem 2 we get a SLag fibration on $M'$. In this case the compactification of the fibers in $M$ is quite different from the case of induced metrics. We will consider one case of this setup in example 1 in Section 2.2. 
\subsection{Examples}

1) {\bf K(N)} (\cite{Gold2}) We will demonstrate a family of examples for the construction in Section 2.1. For those examples we will write down Ricci-flat metrics on $M'$ and we show how fibers of the SLag fibration on $M'$ compactify in $M$. We will also write down some metrics on $M'$, induced from metrics on $M$. We will show that for those metrics the conclusions of Lemma 1 do not apply for non-regular values of the moment map $\mu$ on $D_-$ and thus Lemma 1 is sharp.

Let $N^{2n}$ be a toric Kahler-Einstein manifold with positive scalar curvature $t$. Consider $K(N)$ to be the total space of it's canonical bundle. Then $K(N)$ is a Calabi-Yau manifold with a natural holomorphic $(n+1,0)$-form $\varphi$.
 We refer to \cite{Gold2} for definitions and properties of all relevant structures on $K(N)$. We will however give a Calabi construction of Kahler metrics on $K(N)$, since we will need them in further discussion. Let $r^2 : K(N) \mapsto \mathbb{R}_+$ be the square of the length of elements in $K(N)$. Take any positive function $u$ on $\mathbb{R}_+$ s.t. $u'>0$. Define the metric $\omega_u$ on $K(N)$ as follows: Let $\omega$ be the K-E metric on $N$. The connection on $K(N)$ induces a horizontal distribution for the projection $\pi : K(N) \mapsto N$. We define the horizontal and the vertical distributions to be orthogonal. On the horizontal distribution we define the metric to be $u(r^2)\pi^{\ast}(\omega)$. On the vertical distribution the metric is $t^{-1}u'(r^2)\omega'$. Here $\omega'$ is the induced metric on the linear fibers. The K-E condition ensures that $\omega_u$ is closed. 
If we choose $u(r^2)=(tr^2+l)^{1/n+1}$ for some positive constant $l$, then $\omega_u$ is Ricci-flat. 

Let now $W$ be a 2-plane bundle over $N$, obtained by adding a trivial bundle to $K(N)$. Let $M$ be the projectivization of $W$. The $T^n$-action on $N$ induces an action on $W$ and on $M$. $K(N)$ naturally sits inside of $M$ and $M$ is obtained by adding to $K(N)$ a copy $N_{\infty}$ of $N$ at infinity.
The holomorphic $(n+1,0)$-form $\varphi$ on $K(N)$ becomes a meromorphic $T^n$-invariant section of the canonical bundle of $M$, which defines a divisor $D=D_-=N_{\infty}$. $\varphi$ has singularity of order 2 at $D$ and we get $M'=M-D=K(N)$.

In \cite{Gold2} we considered the SLag fibration on $K(N)$ with respect to the Ricci-flat metric on $K(N)$ (i.e. then $u(r^2)= (tr^2+l)^{1/n+1}$). We have shown that all fibers are asymptotic at infinity to a certain conical fiber $L_0$. Moreover the boundary in $M$ of each fiber is a certain minimal Lagrangian submanifold $L_{\infty}$ of $N_{\infty}$. 

We wish to consider a metric on $K(N)$, induced from a metric on $M$. This can be done by choosing a different function $u$ : Consider a function $w(x)= x^{-2} u'(1/x)$. This is a smooth positive function for $x>0$. Suppose it extends to a smooth positive function $w$ for $x \geq 0$. So in particular $u'(y)$ is asymptotic to $1/y^2$ at infinity and $u'$ integrates on $(0,\infty)$. Thus we have a limit $u_{\infty}>0$ of $u$ at infinity. One can easily show that the metric $\omega_u$ compactifies to a smooth $T^n$-invariant metric $\omega_u$ on $M$.

In \cite{Gold2} we have computed a specific moment map $\mu$ for the $T^n$-action on $N$. Moreover we have shown that the moment map of the action on $K(N)$ was $\mu'=u(r^2) \pi^{-1}(\mu)$. This of course compactifies to a moment map on $M$. Let $\nu$ be a regular value of the moment map $\mu$ on $N$. Then $\nu'=u_{\infty}\nu$ is a regular value of the moment map on $D=N_{\infty}$. Lemma 1 tells us that the level set $D_{\nu'}$ of $\mu'$ on $D$ is a boundary in $M$ of any SLag submanifold of the form $(\mu'=\nu' ~ , ~ Imf=c)$ on $K(N)$.  Let $\nu$ now be a value on the boundary of the moment polytope of $N$. Then we know from $\cite{Gold2}$ that $0$ is in the open part of the moment polytope. Hence by convexity of the polytope we deduce that a multiple $\lambda \mu$ is not in the moment polytope for any $\lambda > 1$. Consider now the value $\nu'= u_{\infty} \nu$ of $\mu'$ on $D$. Then $\nu'$ is not attained by $\mu'$ on $K(N)$. Hence the level set $D_{\nu'}$ does not intersect the closure of any element of our SLag fibration on $K(N)$. Hence the statement of Lemma 1 is sharp.

2) {\bf Toric Varieties}: In section 2.1 we assumed that $H^1(M',\mathbb{R})=0$. We will now demonstrate a class of examples there this condition doesn't hold but we nevertheless have a SLag fibration on $M'$. 

Let $M$ be toric i.e. $T^{n-1} \subset T^n$ and there is a $T^n$-action on $M$.
From Theorem 1 we know that there is a section $\sigma$ of $K(N)$, which is invariant under $T^n$-action.
First we prove the following Lemma, which will also be useful later:
\begin{lem}
Let $N^{2n}$ be a connected complex manifold with a (non-zero) holomorphic $(n,0)$ form $\sigma$ and a Kahler form $\omega$.
Suppose that we have a holomorphic $T^c=(\mathbb{C}^{\ast})^n$-action on $N$ s.t. the action of the real torus $T^n$ is effective, preserves $\sigma$ and Hamiltonian with respect to $\omega$. Then the $T^c$-action is free and $N$ is biholomorphic to $T^c$ under the action. $\sigma$ is non-vanishing and equal to a constant multiple of the form $\bigwedge dz_i / \prod z_i$ on $T^c$. 
Moreover for a choice of $\theta \in \mathbb{R}$ we have that orbits of the $T^n$-action give a SLag fibration on $(N,\omega,e^{i \theta} \sigma)$.
\end{lem}
{\bf Proof}: Let $v_1,\ldots,v_n$ be a basis for the Lie algebra of $T^n$ and let $X_i$ be the flow vector field of $v_i$ on $N$. 
Let $g=i_{X_1}\ldots i_{X_n}\sigma$. As we saw in \cite{Gold2}, $dg=0$, i.e. $g$ is locally constant, and since $N$ is connected $g$ is constant. Since the action is effective, the differential of the action is injective at some point $p$, in which $\sigma(p) \neq 0$. The vectors $X_1(p),\ldots,X_n(p)$ span a Lagrangian plane, so $g=g(p)=g_0 \neq 0$. So $\sigma$ is non-vanishing and the differential of the $T^n$-action is everywhere injective. We can choose
$\theta$ s.t. $e^{ i \theta} g_0$ is real, and we get that orbits of the $T^n$-action are SLag submanifolds of $(N,\omega,e^{i \theta} \sigma)$. 
 
The differential of the $T^c$-action is an isomorphism everywhere. Orbits of the $T^c$-action are open, and since $N$ is connected, there is one orbit. So if $H$ is a stabilizer of some $p \in N$ under $T^c$-action, then $H$ is the stabilizer of $N$ under this action. We wish to prove $H$ is trivial. 
Let $1 \neq h \in H$. Then $h=(z_1,\ldots,z_n)$ with $z_j=a_j + i b_j \neq 0$. Since the $T^n$-action is effective, not all $b_j$ are equal to $0$. The action of $h$ is given by a time 1 flow of the vector field $X=\Sigma a_jX_j + b_jJX_j$. Let $\mu$ be the moment map of the $T^n$-action. Consider a function $\kappa= \Sigma b_j \mu(v_j)$. Then one easily shows that $X(\kappa) >0$, hence the time 1 flow of $X$ cannot return to the same point- a contradiction. Hence $H$ is trivial, the action is free and $T^c$ is biholomorphic to an $N$ under the action. 

In the proof of Lemma 1 we showed that $\sigma$ is invariant under the $T^c$-action since it is invariant under the $T^n$-action. Consider the form $\sigma_0= \bigwedge dz_j / \prod z_j$ on $T^c$. Then both the pullback of $\sigma$ under the action and $\sigma_0$ are $T^c$-invariant. Hence one is a constant multiple of the other and we are done.  Q.E.D.

We return now to the toric manifold $M$. Consider $N=M-D_-$ with a holomorphic $(n,0)$-form $\sigma$ on it. We have a $T^c$-action on $M$, and it leaves $M-D_-$ invariant, hence it induces a $T^c$-action on $N$. From Lemma 2 we deduce that $\sigma$ is non-vanishing, i.e. $D_+ = \emptyset$. Also for a choice of $\theta$ orbits of the $T^n$-action give a SLag fibration on $(M-D_-,\omega,e^{ i \theta} \sigma)$.

We remark that $M'$ coincides with the regular points of the $T$-action. Indeed
points in $M-D_-$ are regular points of $T^n$-action. Also every point $d \in D_-$ is a singular point of the $T$-action (since otherwise the $T^c$-orbit of $d$ is open in $M$, and it cannot be contained in $D_-$).

3) {\bf The Grassmanian G(2,4)}: This is yet another example, for which $H^1(M',\mathbb{R}) \neq 0$ but we can construct a SLag fibration on $M'$.
Consider the Grassmanian $G(2,4)$ of complex 2-planes in $\mathbb{C}^4$, which we identify with a quadric hypersurface $M=(z_1z_2+z_3z_4+z_5z_6=0)$ in $\mathbb{C}P^5$ (see \cite{Wit} and \cite{Hub}). The fourth power $(\gamma^{\ast})^{\otimes 4}$ of the hyperplane bundle on $\mathbb{C}P^5$ restricted to $M$ is the anti-canonical bundle $\overline{K}(M)$. Thus polynomials of degree 4 give rise to holomorphic sections of $\overline{K}(M)$. There is a complex 3-torus $T^c=(\mathbb{C}^{\ast})^3$-action on $M$ given by \[(\lambda_1,\ldots,\lambda_3)(z_1,\ldots,z_6)=(\lambda_1z_1,\lambda_1^{-1}z_2,\ldots,\lambda_3z_5,\lambda_3^{-1}z_6) \]
This action of course contains the action of the real torus $T^3 \subset (\mathbb{C}^{\ast})^3$. We would like to find a homogeneous polynomial $p$ of degree 4 on $\mathbb{C}^6$ s.t. $p$ defines a $T^c$-invariants section of $\overline{K}(M)$. To do that we will carefully set up the $T$-equivariant identifications between various bundles we use. 

First on $\mathbb{C}P^5$ there is a constant bundle $C=\mathbb{C}^6$. It contains a universal bundle $\gamma$ as a sub-bundle. The tangent bundle to $\mathbb{C}P^5$ is isomorphic to $\gamma^{\ast} \otimes (C/\gamma)$. There is a short exact sequence \[\gamma \otimes \gamma^{\ast} \mapsto C \otimes \gamma^{\ast} \mapsto (C/\gamma) \otimes \gamma^{\ast} \]
Taking the canonical bundles of the elements in the sequence we get
$K(\mathbb{C}P^5) \simeq K(C) \otimes \gamma^{\otimes 6}$. Of course the bundle $K(C)$ is trivial, but the isomorphism $ K(\mathbb{C}P^5) \mapsto \gamma^{\otimes 5}$ is not $GL(6,\mathbb{C})$-equivariant, but it is $SL(6,\mathbb{C})$-equivariant. Since our torus $T^c$ is in $SL(6,\mathbb{C})$, we are fine. 

Next consider the quadric $M$, which can be viewed as a zero set of a section $\eta$ of $(\gamma^{\ast})^{\otimes 2}$. The canonical bundle of $M$ is isomorphic to $N \otimes K(\mathbb{C}P^5)$. Here $N$ is the normal bundle to $M$ in $\mathbb{C}P^5$ and the isomorphism is given by $v \otimes \varphi \mapsto i_v \varphi|_{M}$. Also $N$ is isomorphic to $(\gamma^{\ast})^{\otimes 2}$ with an isomorphism given by $v \mapsto \nabla_v \eta$. Since $\eta$ is $T^c$-invariant, this isomorphism is $T^c$-equivariant. So overall we get that $K(M) \simeq N \otimes K(\mathbb{C}P^5) \simeq (\gamma^{\ast})^{\otimes 2} \otimes \gamma^{\otimes 6} \simeq \gamma^{\otimes 4}$. Dually $\overline{K}(M) \simeq (\gamma^{\ast})^{\otimes 4}$ and this isomorphism is $T^c$-equivariant. 

Consider a polynomial $p=(z_1z_2)^2 + (z_3z_4)^2-(z_5z_6)^2$. Then it defines a $T^c$-invariant holomorphic section of $\overline{K}(M)$, and dually it defines a $T^c$-invariant meromorphic section $\sigma$ of $K(M)$, which defines a divisor $D=D_-$, which is the zero set of $p$. Let $M'=M-D$ and pick a $T^3$-invariant Kahler metric $\omega$ on $M'$ s.t. the action is Hamiltonian (e.g one can pick the induced metric from $\mathbb{C}P^5$). We would like to get a SLag fibration for $(M',\omega,\sigma)$. We will see that $H^1(M',\mathbb{R}) \neq 0$, but we can still do that if we replace $\sigma$ by a multiple $e^{i \theta} \sigma$.

Let $e_1,e_2,e_3$ be a basis for Lie algebra of $T^3$ and $X_i$ be corresponding flow vector fields on $M$. Let $\sigma'= i_{X_1} \ldots i_{X_3} \sigma$. Then from \cite{Gold2} we know that $\sigma'$ is a holomorphic $(1,0)$-form on $M'$. Moreover suppose that for some choice of $\theta$ we have $Im(e^{i \theta}\sigma')$ is an exact 1-form on $M'$, i.e. it is equal to $dh$ for some  function $h$. Then we saw in \cite{Gold2} that $(M',\omega,e^{ i \theta}\sigma)$ has a SLag fibration on it with fibers given as level sets of $(\mu,h)$ (here $\mu$ is a moment map for the $T^3$-action on $M'$).

To find this $\theta$ consider the following map $\beta : M' \mapsto \mathbb{C}P^1$ given by $\beta(z_1,\ldots,z_6)=(z_1z_2,z_3z_4)$. We will show that $\beta$ is well-defined and compute it's image. First we claim that on $M'$ we can't have $z_1z_2=0$ or $z_3z_4 = 0$. Indeed suppose $z_1z_2=0$. Then $z_3z_4=-z_5z_6$, hence $p(z_1,\ldots,z_6)=(z_1z_2)^2+(z_3z_4)^2-(z_5z_6)^2=0$ and we are on $D$ and not $M'$.
So we deduce that $\beta$ is well-defined and it's image lies in $C'=\mathbb{C}P^1-((1,0) \cup (0,1))$.

Next we show that $\beta$ is surjective onto $C'$. Indeed $\beta$ has local inverses $\alpha(a,b) = (a,1,b,1,-a-b,1)$ into $M'$ (here $(a,b)$ is in $C'$). $p(a,1,b,1,-a-b,1)=-2ab \neq 0$ and so indeed $(a,1,b,1,-a-b,1) \in M'$. 

The map $\beta$ is $T^c$-invariant.
Take now $(a,b) \in C'$ s.t. $a+b \neq 0$ i.e. $(a,b) \neq (1,-1)$. If $(z_1,\ldots,z_6) \in \beta^{-1}((a,b))$ then all $z_i \neq 0$. One easily deduces that $\beta^{-1}((a,b))$ is in fact the orbit of $\alpha(a,b)$ under $T^c$-action. Moreover because of the local inverse $\alpha$ the differential of $\beta$ is surjective at all points of $\beta^{-1}((a,b))$. So we get a principal fiber bundle $M''=\beta^{-1}(C'-(1,-1))$ over $C'-(1,-1)$ with the fiber being $T^c$.

Now the form $\sigma'$ is $T^c$ -invariant. Moreover the tangent space to the orbits of $T^c$ is in the kernel of $\sigma'$. From this we easily deduce that there is a holomorphic $(1,0)$-form $\sigma''$ on $C'-(1,-1)$ s.t. $\beta^{\ast}(\sigma'')=\sigma'$. We claim that $\sigma''$ compactifies to a holomorphic $(1,0)$-form on $C'$. Indeed we have a local inverse $\alpha$ near $(1,-1)$ and $\alpha^{\ast}(\sigma')=\sigma''$. 

Now $C'$ is $\mathbb{C}^{\ast}$, so $H^1(C',\mathbb{R})=\mathbb{R}$. So one can find an angle $\theta$ s.t. $Im(e^{i \theta}\sigma'')$ is an exact 1-form on $C'$. We claim that also $Im(e^{ i \theta} \sigma')$ is an exact 1-form on $M'$. Indeed take a loop $\gamma$ on $M'$. Then one easily sees that \[ \int_{\gamma}Im(e^{i \theta} \sigma') =  \int_{\beta \circ \gamma} Im(e^{i \theta}\sigma'')=0 \] 
and we are done. 
\section{Calabi-Yau hypersurfaces near a large complex limit Hypersurface}
\subsection{Deformation of SLag fibrations}
In this section we wish to study Calabi-Yau hypersurfaces in an algebraic complexity one space $M$. We will assume that the anti-canonical bundle $\overline{K}$ of $M$ is ample and has a holomorphic section $\eta$, invariant under the $T$-action. Let $D'$ be a smooth part of the divisor $D$ of $\eta$. We will assume that $\eta$ is transversal to $0$ at $D'$. This defines a natural $T^{n-1}$-invariant trivialization $\varphi'$ of the canonical bundle of $D'$. Indeed the canonical bundle of $D'$ is naturally isomorphic to $K \otimes N$. Here $N$ is a normal bundle to $D'$ and the isomorphism is explicitly given by $ \varphi \otimes v \mapsto i_v\varphi|_{D}$. Also the normal bundle $N$ is naturally isomorphic to $\overline{K}$, with isomorphism given by $ v \mapsto \nabla_v \eta$. 

Take a connected component $D_0$ of $D'$.
We assume that the $T^{n-1}$-action on $D_0$ has a finite stabilizer $H$. We have an action of the quotients $T_0=T^{n-1}/H$ and of $T_0^c=T^c/H$ on $D_0$. But $T_0$ is isomorphic to $T^{n-1}$. Also the isomorphism between $T_0$ and $T^{n-1}$ induces a biholomorphic isomorphism between $T_0^c$ and $T^c$. From these isomorphisms we get a holomorphic $T^c$-action on $D_0$, and this action preserves $\varphi'$. Also the corresponding $T^{n-1}$-action is effective and it is Hamiltonian with respect to the metric $\omega$, induced from $M$. From Lemma 2 we deduce that the $T^c$-action is free and $D_0$ is biholomorphic to $T^c$. Also we can choose an angle $\theta_0$ s.t. orbits of the $T^{n-1}$-action give a SLag fibration on $D_0$ for $(\omega, e^{i \theta_0}\varphi')$. We have the following theorem :
\begin{thm}
Let $\eta,D,D',D_0$ be as above. Take $N$ to be any $T$-invariant compact subset of $D_0$ and $U$  a tubular neighbourhood of $N$ in $M$. Then there is a neighbourhood $U_N$ of $\eta$ in $H^0(\overline{K}(M),M)$ s.t. for $\eta_p \in U_N$
there is a SLag fibration of the neighbourhood $U_p$ in $D_p$ with respect to $(\omega_p,e^{i \theta_p}\varphi'_p)$.

Here $D_p$ is the zero set of $\eta_p$, $U_p$ contains $D_p \bigcap U$,
$\omega_p$ is the induced metric on $D_p$ from $M$, $\varphi'_p$ is the natural trivialization of the canonical bundle of $D_p$, given by the section $\eta_p$ and $\eta_p \mapsto \theta_p$ is a smooth function.   
\end{thm}
{\bf Proof:} Choose a compact, $T$-invariant $N \subset D_0$ and pick a compact, $T$-invariant neighbourhood $N'$ of $N$. Let $U'$ be a tubular neighbourhood of $N'$ in $M$.

For any $\eta_p \in H^0(\overline{K}(M),M)$ let $D_p$ be it's zero set. 
The part of $D_p$ in $U'$ has a trivialization $\varphi'_p$ of it's canonical bundle, induced from the section $\eta_p$. Also one can choose a neighbourhood $V$ of $\eta$ in $H^0(\overline{K}(M),M)$ s.t. there is  a map  
$\alpha: N' \times V \mapsto U'$ so that $D_p \bigcap U' = \alpha(N , \eta_p)$ for $\eta_p \in V$. 

We wish to deform the SLag torus fibration on $N'$ to a SLag torus fibration on $D_p \bigcap U$ for a tubular neighbourhood $U \subset U'$ of $N$. 
To set up the deformation theory we consider ${\cal M}$- the connected component of the moduli-space of 
 of $C^{2,\alpha}$ embeddings of $T$ into $N'$, which contains embeddings given by orbits of the $T$-action on $N'$.
This is a Banach manifold, whose local chart at a particular (smooth) embedding are $C^{2,\alpha}$ sections of the normal bundle of the embedded $T$. We consider the space ${\cal M}'= {\cal M} \times V$, which can be thought as a space of embeddings of $L$ into various $D_p$ (with $(f,\eta_p) \mapsto \alpha(f,\eta_p)$ for an embedding $f \in {\cal M}$). 

So let $(f,\eta_p) \in {\cal M}'$ and $\alpha(f,\eta_p)$ be the corresponding embedding of $T$ into $D_p$. Then for a fixed $\eta_p$, various embeddings corresponding to different $f$'s are all isotopic in $D_p$, so they all carry the same isotopy class $L_p$ in $D_p$.
The cohomology class $[\omega]$ of $\omega$ restricted to $D_p$ is the cohomology class induced from $\omega$ on $M$. Also the image of $\alpha(f,\eta_p)$ is isotopic in $M$ to an orbit of the $T$-action on $N'$ and $\omega$ restricts to $0$ on the orbits. So
$[\omega]|_{L_p}=0$. Also on $N$ $\varphi'$ evaluated on orbits was not zero.
By continuity $\varphi'_p(L_p) \neq 0$. So we can choose $e^{i \theta_p}$ s.t.
\[Im(e^{i \theta_p}\varphi'_p(L_p))=0\] This $e^{i \theta_p}$ is defined up to a sign, and we can choose it to be a smooth function of $\eta_p$.

We want to consider the subset $SLag \subset {\cal M}'$ of those embeddings $(f,\eta_p)$ s.t. the forms $\omega$ and $Im(e^{i \theta_p}\varphi'_p)$ restrict to $0$ on $\alpha(f,\eta_p)$. Let $\pi: {\cal M}' \mapsto V$ be the projection onto the second factor. We also have a space $S_0$ of $T$-orbits in $N'$, and we will think of $S_0 \subset \pi^{-1}(\eta) \bigcap SLag$. We have the following:
\begin{lem}
The space $SLag$ near $S_0$ is a smooth manifold $SLag_0$ of dimension $n-1 +2dim(H^0(\overline{K}(M),M))$. Moreover the differential of the projection $\pi$, restricted to $SLag_0$ is surjective.
\end{lem}
We claim that the statement of the Theorem follows from this Lemma. Indeed $SLag_0$ will be a fibration over $V'$ for sufficiently small neighbourhood $V'$ of $\eta$. If we choose any metric on $SLag_0$ (for instance the one induced from the product metric on ${\cal M}$ and $V'$), then we have a horizontal distribution for the projection $\pi$ on $SLag_0$. For $\eta_p \in V'$ we can look at a line $t \mapsto \eta + t(\eta_p-\eta)$. This line induces a flow $\rho_p(t)$ on $SLag_0$, s.t $\rho_p(1)$ sends $\pi^{-1}(\eta)$ to $\pi^{-1}(\eta_p)$. It is an easy exercise in differential topology that one can choose a smaller neighbourhood $U_N$ of $\eta$ in $H^0(\overline{K}(M),M)$ s.t for $\eta_p \in U_N$ the image of $S_0$ under $\rho_p(1)$ gives a SLag fibration, which covers the neighbourhood $U \bigcap D_p$.   

Now we prove the {\bf Lemma}: Since $S_0$ is compact, it is obviously enough to prove our claim near a point $L \in S_0$ (here we think of $L$ as an orbit of the $T$-action in $N'$). A neighbourhood $Y$ of $L$ in ${\cal M}$ can be thought as a small ball in the space of $C^{2,\alpha}$ normal vector fields to $L$. There is a Banach vector bundle over ${\cal M}'$, given by the direct sum of exact $C^{1,\alpha}$ 2-forms with exact $C^{1,\alpha}$ $(n-1)$-forms on $T^{n-1}$. There is a section $\sigma$ of this bundle over $Y \times V$, given by $\sigma(f,\eta_p)= (\alpha(f,\eta_p)^{\ast}(\omega),\alpha(f,\eta_p)^{\ast}(Im(e^{i \theta_p}\varphi'_p)))$. The space $SLag$ precisely corresponds to the zero set of this section.

We have shown in \cite{Gold1} that one can slightly generalize the argument of \cite{Mc} to prove that the differential of this $\sigma$ is already surjective then restricted to the tangent space $T_L{\cal M} \subset T_L{\cal M}'$ and the kernel of $d\sigma$ in $T_L{\cal M}$ is of dimension $n-1$ (in our case the kernel
in $T{\cal M}$ corresponds to vectors fields, which realize deformations of the orbits). So one deduces that $SLag$ is smooth of dimension $n-1+2dim(H^0(\overline{K}(M),M))$ near $L$ and the differential of $\pi$ restricted to $SLag$ is surjective. Q.E.D.  

{\bf Remark 2}: Suppose $D'$ has several connected components. For each connected component we have a SLag fibration on a corresponding neighbourhood of $D_p$. But those fibrations are with respect to holomorphic volume forms
$e^{i \theta_p} \varphi'_p$ on $D_p$, there $\theta_p$ are potentially different for various connected components of $D'$. Suppose $D_p$ is smooth. Then $\varphi'_p$ is a holomorphic volume form on $D_p$.
 It is clear from the proof that $\theta_p$ will coincide if the corresponding isotopy classes $L_p$ coming from different connected components of $D'$ define the same homology class in $D_p$. We will see two examples (3.2 and 3.3), in which this holds. Thus in particular we will produce SLag fibrations on subsets of $D_p$, whose complements have arbitrary small volume in $D_p$.

We would like to point out that Theorem 3 generalizes to complete intersections on a Kahler manifold $M$ with a $T$-action of complexity larger then 1. Indeed let $(M^{2n},\omega)$ be a Kahler manifold. Suppose that the anticanonical bundle $\overline{K}(M)$ is a tensor product 
\begin{equation}
\overline{K}(M) \simeq L_1 \otimes \ldots \otimes L_d
\end{equation}
Suppose we have a $T^{n-d}$-action on $M$ s.t $T^{n-d}$ also acts on the total space of each $L_i$. Moreover we assume that the action is equivariant with respect to the isomorphism (1). Suppose we have sections $\eta_i$ of $L_i$, which are invariant under the $T$-action. Let $D$ be their common zero set and $D'$ be the smooth part of $D$. We assume that $(\eta_1,\ldots,\eta_d)$ is transversal to $0$ along $D'$. Thus we get a trivialization $\varphi'$ of the canonical bundle of $D'$ and $\varphi'$ is invariant under the $T$ and the $T^c$-actions (here $T^c$ is a complex torus of dimension $n-d$). 

Let $D_0$ be a connected component of $D'$. Assume that the $T$-action on $D_0$ has a finite stabilizer $H$. We have, as before, the actions of the quotients $T/H$ and $T^c/H$ on $D_0$. The $T^c/H$ action will be free, $D_0$ is biholomorphic to $T^c$ and we can choose an angle $\theta_0$ s.t. orbits of the $T$-action will be SLag submanifolds of $(D_0,\omega,e^{i \theta_0}\varphi')$. We will also have the following result, which generalizes Theorem 3:
\begin{cor}
Let $\eta=(\eta_1,\ldots,\eta_d),D,D',D_0$ be as above. Let $N$ be a compact, $T$-invariant subset of $D_0$ and $U$ be it's tubular neighbourhood in $M$. Then we can choose a neighbourhood $U_N$ of $\eta$ in $H^0(L_1,M) \times \ldots \times H^0(L_d,M)$ s.t for any $\eta^p=(\eta_1^p,\ldots,\eta_d^p) \in U_N$ we have a SLag torus fibration on a neighbourhood $U_p$ of $D_p$ with respect to $(\omega,e^{i \theta_p}\varphi'_p)$.

Here $D_p$ is the common zero set of $(\eta_1^p,\ldots,\eta_d^p)$, $U_p$ contains $U \bigcap D_p$, $\omega$ is the induced metric on $D_p$, $\varphi'_p$ is the trivialization of the canonical bundle of $D_p$, given by the section $\eta_p$ and $\eta_p \mapsto \theta_p$ is smooth.
\end{cor}
The proof of the Corollary is analogous to the proof of Theorem 3.

In Sections 3.2-3.4 we consider some examples there Theorem 3 and Corollary 1 apply.
\subsection{Fermat type quintics in $\mathbb{C}P^4$}
Consider $\mathbb{C}P^n$ with a $T^c=(\mathbb{C}^{\ast})^{n-1}$-action given by \[(c_1,\ldots,c_{n-1}) (z_1,\ldots,z_{n+1})= (c_1 z_1,\ldots,c_{n-1} z_{n-1},\prod c_j^{-1} z_n,z_{n+1}) \]
This action of course contains a $T^{n-1} \subset T^c$-action.
Consider the following polynomial $f_{\lambda}= \prod z_j - \lambda (z_{n+1})^{n+1}$. Then for every $\lambda$ this polynomial defines a section of the anti-canonical bundle $\overline{K}(\mathbb{C}P^n)$, which is invariant under the $T^c$-action. There are essentially 2 distinct cases $\lambda=0$ or $\lambda=1$.

${\bf \lambda= 0}$: $D$ has $n+1$ connected components $D_j$, given by $D_j= (z_j=0)$. Let $D'_j$ be the smooth part of $D_j$ in $D$, i.e. $D'_j=(z_j=0 ~ , ~z_i \neq 0 ~ for ~ i \neq j)$.
 Consider the following 1-parameter family $f^t=f_0+ t\Sigma z_j^{n+1}$. Then $n=4$ the zero sets $D_t$ are Fermat type quintics. From each component $D'_j$ we get a SLag fibration on a corresponding part of $D_t$, but the fibrations are SLag for potentially different holomorphic volume forms on $D_t$ for different components $D'_j$ (see Remark 2). We will show that the isotopy classes $L_t^j$, coming from components $D'_j$, give the same homology class in $D_t$. Thus SLag fibrations will be with respect to one holomorphic volume form $\varphi_t$ on $D_t$.

W.l.o.g it is enough to show that $L_t^1$ and $L_t^2$ define the same homology class in $D_t$. To show this consider a different $T^{n-1}$-action on $\mathbb{C}P^n$, given by \[ (e^{ i \theta_1},\ldots,e^{ i \theta_{n-1}})(z_1,\ldots,z_{n+1})=(z_1,z_2,e^{i \theta_1}z_3,\ldots, e^{i \theta_{n-1}}z_{n+1}) \]
This action has a moment map $\mu=(\mu_1,\ldots,\mu_{n-1})$ with \[ \mu_j= \frac{|z_{j+2}|^2}{\Sigma |z_i|^2} \]

Consider the level set $\Sigma=\mu^{-1}(1/2n,\ldots,1/2n)$ restricted to $D$. It's intersection with each of $D'_1$ and $D'_2$ it will be 1 smooth orbit $L_i$ in the interior ($i=1,2$).  Also $\Sigma$ will not intersect $D_j$ for $j>2$.  Also $L_i$ are easily shown to coincide with the orbits of the original $T^{n-1}$-action on $D$.

Consider now $D_t$ and $\mu|_{D_t}$. Then the level set $\Sigma_t=\mu^{-1}(1/2n,\ldots,1/2n)$ on $D_t$ will also have 2 smooth components. Also clearly those components will be in the isotopy class of $L_t^1$ and $L_t^2$. Now the level set of any map to a non-compact manifold carries a trivial homology class. So $L_t^1$ and $L_t^2$ (with corresponding orientations) will have the same homology class in $D_t$ and we are done.

{\bf Remark 3}: We have constructed Special Lagrangian fibrations on large parts of quintics near the large complex limit quintic. If one forgets about the Special condition and studies Lagrangian fibrations, then one can say more. In fact W.-D. Ruan has constructed in \cite{Ruan} and \cite{Ruan2} Lagrangian tori fibrations on general quintics and also constructed symplectic mirrors to those fibrations.

${\bf \lambda=1}$: The singular points of $D$ are points there $z_{n+1}=0$ and $z_1 \cdots z_n=0$. The set of smooth points $D'$ will have 2 connected components: $D'_1=(z_{n+1}=0 ~ , ~ z_j \neq 0 ~ for ~ j<n+1)$ and $D'_2=(z_j \neq 0)$. $D'_1$ will be the orbit $(1,\ldots,1,0)$ under the $T^c$-action; $D'_2$ will be the orbit of $(1,\ldots,1,1)$ under $T^c$-action. 

We have a moment map $\mu=(\mu_1,\ldots,\mu_{n-1})$ for the $T^{n-1}$-action on $\mathbb{C}P^n$ with $\mu_j= \frac{|z_j|^2-|z_n|^2}{\Sigma|z_j|^2}$. One easily sees that the preimage $\mu^{-1}(0,\ldots,0)$ doesn't intersect the singular points of $D$. Also it intersects each component $D'_j$ in one smooth orbit of the $T^{n-1}$-action. So we use the same trick as before to show that the isotopy classes $L_t^1$ and $L_t^2$ define the same homology class on $D_t$. Thus
from Theorem 3 we get SLag fibrations on large parts of Calabi-Yau hypersurfaces near the complex limit hypersurface $D$.
\subsection{Calabi-Yau Hypersurfaces in the Grassmanian $G(2,4)$}
Consider the Grassmanian $M=G(2,4)$ of 2-planes in $\mathbb{C}^4$, which we continue to identify with a quadric hypersurface $z_1z_2+z_3z_4+z_5z_6=0$ in $\mathbb{C}P^5$. We have a $T^3$ and a $T^c=(\mathbb{C}^{\ast})^3$-action on $M$ as in Section 2.2. The polynomial $f=(z_1z_2)^2+(z_3z_4)^2+(z_5z_6)^2$ defines a $T^c$-invariant section of the anti-canonical bundle $\overline{K}(M)$. 

The singular part of the zero set $D$ of $f$ corresponds to points $z_{2j-1}z_{2j}=0$ for $j=1,2,3$. Let $D'$ be the smooth part of $D$ and $(z_1,\ldots,z_6) \in D'$. We can assume w.l.o.g. that $z_1z_2=1$. Let $z_3z_4=a$. Then $z_5z_6=-a-1$ and $a$ satisfies $a^2+a+1=0$.
We get 2 values $a_i= \frac{-1 \pm \sqrt{3}i}{2}$. For each $a_i$ we have a point $b_i=(1,1,1,a_i,1,-1-a_i)$ on $D'$. The orbits of $T^c$ through $b_i$ give 2 connected components $D'_1$ and $D'_2$ of $D'$ (explicitly $D'_i$ are given by $z_3z_4/z_1z_2=a_i$).

The $T^3$-action has a moment map $\mu=(\mu_1,\mu_2,\mu_3)$ with $\mu_i=\frac{|z_{2i-1}|^2-|z_{2i}|^2}{\Sigma|z_j|^2}$. The preimage $\mu^{-1}(0,0,0)$ doesn't intersect the singular part of $D$. Also it intersects  each of $D'_j$ at one smooth orbit. We use the same argument as in Section 3.2 to prove that the isotopy classes $L_p^1$ and $L_p^2$ on $D_p$ give the same homology class. Thus we can use Theorem 3 to obtain SLag fibrations on large parts of Calabi-Yau hypersurfaces $D_p$, which are zero sets of polynomials $f_p$ of degree 4 for $f_p$ sufficiently close to $f$.

{\bf Remark 4}: In \cite{Gold1} we gave an example of a SLag submanifold on a Calabi-Yau hypersurface in $G(2,4)$. 
\subsection{Complete intersection of two degree 3 hypersurfaces in $\mathbb{C}P^5$}
In this section we want to illustrate an application of Corollary 1. Let $M$ be $\mathbb{C}P^5$. We decompose it's anticanonical bundle as 
\begin{equation}
\overline{K}(M) \simeq L_1 \otimes L_2
\end{equation}
Here $L_1=L_2=(\gamma^{\ast})^{\otimes 3}$. We have a $T^3$-action on $M$, given by \[(e^{i \theta_1},\ldots,e^{i \theta_3})(z_1,\ldots,z_6)=(e^{i \theta_1}z_1,e^{i \theta_2}z_2,e^{i \theta_3}z_3,e^{-i(\theta_1+\theta_2)}z_4,e^{i(\theta_1+\theta_2-\theta_3)}z_5,e^{-i(\theta_1+\theta_2)}z_6)\]
Then $T$ acts on $\overline{K}(M),L_1,L_2$. Since the linear action of $T$ on $\mathbb{C}^6$ is in $SL(6,\mathbb{C})$, the $T$-action is equivariant with respect to the isomorphism (2) (see example (3) in Section 2.2).
We have 4 monomials $g_1=z_1z_2z_4 ~ , ~ g_2= z_3z_4z_5 ~ , ~ g_3=z_1z_2z_6 ~ , ~ g_4=z_3z_5z_6$, which can be viewed as $T$-invariant sections of $L_1$ and $L_2$. We pick $\eta_1$ and $\eta_2$ to be some of their linear combinations s.t. conditions of Corollary 1 hold. Thus we get SLag fibrations on large parts of complete intersections of 2 hypersurfaces of degree 3 in $\mathbb{C}P^5$ near $(\eta_1,\eta_2)$.

Massachusetts Institute of Technology

E-Mail : egold@math.mit.edu

\end{document}